%Plain TeX file
%
\magnification\magstep1
\baselineskip15pt
\hfuzz5pt

\newcount\relFnno
\def\ref#1{\expandafter\edef\csname#1\endcsname}

% \newread\AUX\immediate\openin\AUX=\jobname.aux
% \ifeof\AUX\immediate\write16{\jobname.aux gibt es nicht!}\else
% \input \jobname.aux
% \fi\immediate\closein\AUX

\relFnno 1
\ref {Intro}{1}
\ref {Fn1}{\the \relFnno }
\advance \relFnno 1
\relFnno 1
\ref {NullCone}{2}
\ref {TA}{Theorem\penalty 10000\ 2.1}
\relFnno 1
\ref {CZ}{Corollary\penalty 10000\ 2.2}
\relFnno 1
\ref {degrees}{3}
\ref {TB}{Lemma\penalty 10000\ 3.1}
\ref {TC}{Lemma\penalty 10000\ 3.2}
\relFnno 1
\ref {TD}{Theorem\penalty 10000\ 3.3}
\ref {CA}{Corollary\penalty 10000\ 3.4}
\ref {CB}{Corollary\penalty 10000\ 3.5}
\ref {TE}{Theorem\penalty 10000\ 3.6}
\ref {CD}{Corollary\penalty 10000\ 3.7}
\relFnno 1
\ref {TF}{Theorem\penalty 10000\ 3.8}
\ref {DomHeg}{Theorem\penalty 10000\ 3.9}
\ref {DomHegSez}{Theorem\penalty 10000\ 3.10}
\ref {etabeta}{Theorem\penalty 10000\ 3.11}
\relFnno 1
\ref {Functoriality}{4}
\ref {faith}{Theorem\penalty 10000\ 4.1}
\ref {E5}{$(4.2)$}
\relFnno 1
\ref {modular}{Theorem\penalty 10000\ 4.3}
\ref {vector}{Proposition\penalty 10000\ 4.4}
\ref {generalk}{Theorem\penalty 10000\ 4.5}
\ref {E6}{$(4.4)$}
\relFnno 1
\ref {E8}{$(4.5)$}
\ref {E7}{$(4.6)$}
\relFnno 1
\ref {Polarization}{5}
\ref {Polar}{Theorem\penalty 10000\ 5.1}
\ref {E1}{$(5.3)$}
\relFnno 1
\ref {E2}{$(5.7)$}
\relFnno 1
\ref {E4}{$(5.10)$}
\ref {E3}{$(5.12)$}
\ref {WeylTheorem}{6}
\ref {Extend}{Theorem\penalty 10000\ 6.1}
\ref {Fn2}{\the \relFnno }
\advance \relFnno 1
\relFnno 1
\ref {E9}{$(6.1)$}
\relFnno 1
\ref {UniCrit}{Theorem\penalty 10000\ 6.3}
\ref {regrep1}{Corollary\penalty 10000\ 6.5}
\relFnno 1
\ref {improvement}{7}
\ref {onedim}{Theorem\penalty 10000\ 7.1}
\ref {wUniCrit}{Theorem\penalty 10000\ 7.3}
\relFnno 1
\ref {regu}{Corollary\penalty 10000\ 7.4}
\ref {regrep2}{Corollary\penalty 10000\ 7.5}
\ref {References}{8}
\ref {DH}{[DH]}
\ref {Fl}{[Fl1]}
\ref {Fleisch}{[Fl2]}
\ref {Fo}{[Fo]}
\ref {Noe}{[N]}
\ref {Ri}{[R]}
\ref {Sch}{[Sch]}
\ref {Se}{[Se]}
\ref {Sm}{[Sm]}
\ref {Weyl}{[We]}
\ref {Spaltenbreite}{24.1667pt}
\relFnno 1

% \input makros

% Wirkungsweise:

% ...bearbeiteter Text...
% \ignore                   <--- muss am Ende einer Zeile stehen
% ...ignorierter Text...    <--- keine Zeile darf mit ';' beginnen
% ;;                        <--- muss allein auf einer Zeile stehen
% ...bearbeiteter Text...

\def\ignore{\bgroup
\catcode`\;=0\catcode`\^^I=14\catcode`\^^J=14\catcode`\^^M=14
\catcode`\ =14\catcode`\!=14\catcode`\"=14\catcode`\#=14\catcode`\$=14
\catcode`\&=14\catcode`\'=14\catcode`\(=14\catcode`\)=14\catcode`\*=14
\catcode`+=14\catcode`\,=14\catcode`\-=14\catcode`\.=14\catcode`\/=14
\catcode`\0=14\catcode`\1=14\catcode`\2=14\catcode`\3=14\catcode`\4=14
\catcode`\5=14\catcode`\6=14\catcode`\7=14\catcode`\8=14\catcode`\9=14
\catcode`\:=14\catcode`\<=14\catcode`\==14\catcode`\>=14\catcode`\?=14
\catcode`\@=14\catcode`\A=14\catcode`\B=14\catcode`\C=14\catcode`\D=14
\catcode`\E=14\catcode`\F=14\catcode`\G=14\catcode`\H=14\catcode`\I=14
\catcode`\J=14\catcode`\K=14\catcode`\L=14\catcode`\M=14\catcode`\N=14
\catcode`\O=14\catcode`\P=14\catcode`\Q=14\catcode`\R=14\catcode`\S=14
\catcode`\T=14\catcode`\U=14\catcode`\V=14\catcode`\W=14\catcode`\X=14
\catcode`\Y=14\catcode`\Z=14\catcode`\[=14\catcode`\\=14\catcode`\]=14
\catcode`\^=14\catcode`\_=14\catcode`\`=14\catcode`\a=14\catcode`\b=14
\catcode`\c=14\catcode`\d=14\catcode`\e=14\catcode`\f=14\catcode`\g=14
\catcode`\h=14\catcode`\i=14\catcode`\j=14\catcode`\k=14\catcode`\l=14
\catcode`\m=14\catcode`\n=14\catcode`\o=14\catcode`\p=14\catcode`\q=14
\catcode`\r=14\catcode`\s=14\catcode`\t=14\catcode`\u=14\catcode`\v=14
\catcode`\w=14\catcode`\x=14\catcode`\y=14\catcode`\z=14\catcode`\{=14
\catcode`\|=14\catcode`\}=14\catcode`\~=14\catcode`\^^?=14
\Ignoriere}
\def\Ignoriere#1\;{\egroup}

\def\today{\number\day.~\ifcase\month\or
  Januar\or Februar\or M{\"a}rz\or April\or Mai\or Juni\or
  Juli\or August\or September\or Oktober\or November\or Dezember\fi
  \space\number\year}
\font\sevenex=cmex7
%\font\sevenex=cmex10 scaled 700
\scriptfont3=\sevenex
%\font\fiveex=cmex7 scaled 714
\font\fiveex=cmex10 scaled 500
\scriptscriptfont3=\fiveex
\def\A{{\bf A}}
\def\G{{\bf G}}
\def\P{{\bf P}}

\def\phi{\varphi}
\def\epsilon{\varepsilon}
\def\theta{\vartheta}
%\font\lams=lams3
%\def\uinto{\lower1.7pt\hbox{%
%\vbox{\offinterlineskip
%\hbox{\lams\char"7A}%
%\hbox{\vbox to 7.5pt{\leaders\vrule width0.2pt\vfill}%
%\kern-.3pt\hbox{\lams\char"76}}}}}
\def\uauf{\lower1.7pt\hbox to 3pt{%
\vbox{\offinterlineskip
\hbox{\vbox to 8.5pt{\leaders\vrule width0.2pt\vfill}%
\kern-.3pt\hbox{\lams\char"76}\kern-0.3pt%
$\raise1pt\hbox{\lams\char"76}$}}\hfil}}

\def\title#1{\par
{\baselineskip1.5\baselineskip\rightskip0pt plus 5truecm
\leavevmode\vskip0truecm\noindent\font\BF=cmbx10 scaled \magstep2\BF #1\par}
\vskip1truecm
\leftline{\font\CSC=cmcsc10
{\CSC Friedrich Knop}\footnote{Supported by an NSA-grant.}}
\leftline{Department of Mathematics, Rutgers University, New Brunswick NJ
08903, USA}
\leftline{knop@math.rutgers.edu}
\vskip1truecm
\par}

%%%%%%%%%%%%%%%%%%
% Makros f"ur Querverweise:
\def\cite#1{\expandafter\ifx\csname#1\endcsname\relax
{\bf?}\immediate\write16{#1 ist nicht definiert!}\else\csname#1\endcsname\fi}
\def\expandwrite#1#2{\edef\next{\write#1{#2}}\next}
\def\neverexpand{\noexpand\noexpand\noexpand}
\def\strip#1\ {}
\def\ncite#1{\expandafter\ifx\csname#1\endcsname\relax
{\bf?}\immediate\write16{#1 ist nicht definiert!}\else
\expandafter\expandafter\expandafter\strip\csname#1\endcsname\fi}
\newwrite\AUX
\immediate\openout\AUX=\jobname.aux
%%%%%%%%%%%%%%%%%%%%%
% Definition von \eightpoint:
\font\eightrm=cmr8\font\sixrm=cmr6
\font\eighti=cmmi8
\font\eightit=cmti8
\font\eightbf=cmbx8
\font\eightcsc=cmcsc10 scaled 833
\def\eightpoint{%
\textfont0=\eightrm\scriptfont0=\sixrm\def\rm{\fam0\eightrm}%
\textfont1=\eighti
\textfont\bffam=\eightbf\def\bf{\fam\bffam\eightbf}%
\textfont\itfam=\eightit\def\it{\fam\itfam\eightit}%
\def\csc{\eightcsc}%
\setbox\strutbox=\hbox{\vrule height7pt depth2pt width0pt}%
\normalbaselineskip=0,8\normalbaselineskip\normalbaselines\rm}
%%%%%%%%%%%%%%%%%%
% Fussnotenmakros
\newcount\absFnno\absFnno1
\write\AUX{\relFnno1}
\newif\ifMARKE\MARKEtrue
{\catcode`\@=11
\gdef\footnote{\ifMARKE\edef\@sf{\spacefactor\the\spacefactor}\/%
$^{\cite{Fn\the\absFnno}}$\@sf\fi
\MARKEtrue
\insert\footins\bgroup\eightpoint
\interlinepenalty100\let\par=\endgraf
\leftskip=0pt\rightskip=0pt
\splittopskip=10pt plus 1pt minus 1pt \floatingpenalty=20000\smallskip
\item{$^{\cite{Fn\the\absFnno}}$}%
\expandwrite\AUX{\neverexpand\ref{Fn\the\absFnno}{\neverexpand\the\relFnno}}%
\global\advance\absFnno1\write\AUX{\advance\relFnno1}%
\bgroup\strut\aftergroup\@foot\let\next}}
\skip\footins=12pt plus 2pt minus 4pt
\dimen\footins=30pc
\output={\plainoutput\immediate\write\AUX{\relFnno1}}
%%%%%%%%%%%%%%%%%%%%%
\newcount\Abschnitt\Abschnitt0
\def\beginsection#1. #2 \par{\advance\Abschnitt1%
\vskip0pt plus.10\vsize\penalty-250
\vskip0pt plus-.10\vsize\bigskip\vskip\parskip
\edef\TEST{\number\Abschnitt}
\expandafter\ifx\csname#1\endcsname\TEST\relax\else
\immediate\write16{#1 hat sich geaendert!}\fi
\expandwrite\AUX{\neverexpand\ref{#1}{\TEST}}
\leftline{\marginnote{#1}\bf\number\Abschnitt. \ignorespaces#2}%
\nobreak\smallskip\noindent\SATZ1\GNo0}
%%%%%%%%%%%%%%%%%%
\def\Proof:{\par\noindent{\it Proof:}}
\def\Remark:{\ifdim\lastskip<\medskipamount\removelastskip\medskip\fi
\noindent{\bf Remark:}}
\def\Remarks:{\ifdim\lastskip<\medskipamount\removelastskip\medskip\fi
\noindent{\bf Remarks:}}
\def\Definition:{\ifdim\lastskip<\medskipamount\removelastskip\medskip\fi
\noindent{\bf Definition:}}
\def\Example:{\ifdim\lastskip<\medskipamount\removelastskip\medskip\fi
\noindent{\bf Example:}}
\def\Examples:{\ifdim\lastskip<\medskipamount\removelastskip\medskip\fi
\noindent{\bf Examples:}}
%%%%%%%%%%%%%%%%
\newif\ifmarginalnotes\marginalnotesfalse
\newif\ifmarginalwarnings\marginalwarningstrue

\def\marginnote#1{\ifmarginalnotes\hbox to 0pt{\eightpoint\hss #1\ }\fi}

\def\strutdepth{\dp\strutbox}
\def\Randbem#1#2{\ifmarginalwarnings
{#1}\strut
\setbox0=\vtop{\eightpoint
\rightskip=0pt plus 6mm\hfuzz=3pt\hsize=16mm\noindent\leavevmode#2}%
\vadjust{\kern-\strutdepth
\vtop to \strutdepth{\kern-\ht0
\hbox to \hsize{\kern-16mm\kern-6pt\box0\kern6pt\hfill}\vss}}\fi}

\def\Zitat!{\Randbem{\bf?}{\bf Zitat}}

\newcount\SATZ\SATZ1
\def\proclaim #1. #2\par{\ifdim\lastskip<\medskipamount\removelastskip
\medskip\fi
\noindent{\bf#1.\ }{\it#2}\Par
\ifdim\lastskip<\medskipamount\removelastskip\goodbreak\medskip\fi}
\def\Aussage#1{%
\expandafter\def\csname#1\endcsname##1.{\ifx?##1?\relax\else
\edef\TEST{#1\penalty10000\ \number\Abschnitt.\number\SATZ}
\expandafter\ifx\csname##1\endcsname\TEST\relax\else
\immediate\write16{##1 hat sich geaendert!}\fi
\expandwrite\AUX{\neverexpand\ref{##1}{\TEST}}\fi
\proclaim {\marginnote{##1}\number\Abschnitt.\number\SATZ. #1\global\advance\SATZ1}.}}
\Aussage{Theorem}
\Aussage{Proposition}
\Aussage{Corollary}
\Aussage{Lemma}
%%%%%%%%%%%%%%%%
\font\la=lasy10
\def\strich{\hbox{$\vcenter{\hbox
to 1pt{\leaders\hrule height -0,2pt depth 0,6pt\hfil}}$}}
\def\dashedrightarrow{\hbox{%
\hbox to 0,5cm{\leaders\hbox to 2pt{\hfil\strich\hfil}\hfil}%
\kern-2pt\hbox{\la\char\string"29}}}

\def\Bindestrich{\penalty10000-\hskip0pt}
\let\_=\Bindestrich
\def\.{{\sfcode`.=1000.}}
%%%%%%%%%%%%%%%%%%%%%%%%%%%%%%%%%%%%

\def\Par{\par}
\def\:={\mathrel{\raise0,9pt\hbox{.}\kern-2,77779pt
\raise3pt\hbox{.}\kern-2,5pt=}}
\def\=:{\mathrel{=\kern-2,5pt\raise0,9pt\hbox{.}\kern-2,77779pt
\raise3pt\hbox{.}}} 
\def\into{\hookrightarrow}
\def\pfeil{\rightarrow}

\def\Pfeil{\longrightarrow}

\def\Pf#1{\buildrel#1\over\longrightarrow}

\def\Ugleich{\hbox{$\cup$\kern.5pt\vrule depth -0.5pt}}
\def\|#1|{\mathop{\rm#1}\nolimits}
\def\<{\langle}
\def\>{\rangle}
\let\Times=\times
\def\times{\mathop{\Times}}
\let\Otimes=\otimes
\def\otimes{\mathop{\Otimes}}
%%%%%%%%%%%%%%%%%%%%%%%%%%%%%%%%%
%Laden von Fonts:
\catcode`\@=11
\def\hex#1{\ifcase#1 0\or1\or2\or3\or4\or5\or6\or7\or8\or9\or A\or B\or
C\or D\or E\or F\else\message{Warnung: Setze hex#1=0}0\fi}
\def\fontdef#1:#2,#3,#4.{%
\alloc@8\fam\chardef\sixt@@n\FAM
\ifx!#2!\else\expandafter\font\csname text#1\endcsname=#2
\textfont\the\FAM=\csname text#1\endcsname\fi
\ifx!#3!\else\expandafter\font\csname script#1\endcsname=#3
\scriptfont\the\FAM=\csname script#1\endcsname\fi
\ifx!#4!\else\expandafter\font\csname scriptscript#1\endcsname=#4
\scriptscriptfont\the\FAM=\csname scriptscript#1\endcsname\fi
\expandafter\edef\csname #1\endcsname{\fam\the\FAM\csname text#1\endcsname}
\expandafter\edef\csname hex#1fam\endcsname{\hex\FAM}}
\catcode`\@=12 

%%%%%%%%%%%%%%%%%%%%%%%%%%%%%%%%%
\fontdef Ss:cmss10,,.
\fontdef Fr:eufm10,eufm7,eufm5.

                        %Hier aufpassen!!!

\def\fm{{\Fr m}}

\fontdef bbb:msbm10,msbm7,msbm5.
\fontdef mbf:cmmib10,cmmib7,.
\fontdef msa:msam10,msam7,msam5.
\def\CC{{\bbb C}}
\def\FF{{\bbb F}}

\def\NN{{\bbb N}}
\def\QQ{{\bbb Q}}

\def\ZZ{{\bbb Z}}

\def\cP{{\cal P}}

\mathchardef\leer=\string"0\hexbbbfam3F
\mathchardef\subsetneq=\string"3\hexbbbfam24
\mathchardef\semidir=\string"2\hexbbbfam6E
\mathchardef\dirsemi=\string"2\hexbbbfam6F
\mathchardef\haken=\string"2\hexmsafam78
\mathchardef\auf=\string"3\hexmsafam10
\let\OL=\overline
\def\overline#1{{\hskip1pt\OL{\hskip-1pt#1\hskip-.3pt}\hskip.3pt}}
\def\aq{{\overline{a}}}\def\Aq{{\overline{A}}}
\def\Bq{{\overline{B}}}

%<--                    Aufpassen  

%
%%%%%%%%%%%%
\newdimen\Parindent
\Parindent=\parindent

%%%%%%%%%%%%
% Displayroutine

\abovedisplayskip 9.0pt plus 3.0pt minus 3.0pt
\belowdisplayskip 9.0pt plus 3.0pt minus 3.0pt
\newdimen\Grenze\Grenze2\Parindent\advance\Grenze1em
\newdimen\Breite
\newbox\DpBox
\def\NewDisplay#1$${\Breite\hsize\advance\Breite-\hangindent
\setbox\DpBox=\hbox{\hskip2\Parindent$\displaystyle{#1}$}%
\ifnum\predisplaysize<\Grenze\abovedisplayskip\abovedisplayshortskip
\belowdisplayskip\belowdisplayshortskip\fi
\global\futurelet\nexttok\WEITER}
\def\WEITER{\ifx\nexttok\qed\expandafter\leftQEDdisplay
\else\leftdisplay\fi}
\def\leftdisplay{\hskip-\hangindent\leftline{\box\DpBox}$$}
\def\leftQEDdisplay{\hskip-\hangindent
\line{\copy\DpBox\hfill\lower\dp\DpBox\copy\QEDbox}%
\belowdisplayskip0pt$$\bigskip\let\nexttok=}
\everydisplay{\NewDisplay}
%%%%%%%%%%%
\newcount\GNo\GNo=0
\newcount\maxEqNo\maxEqNo=0
\def\eqno#1{%
\global\advance\GNo1%
\edef\FTEST{$(\number\Abschnitt.\number\GNo)$}%
\ifx?#1?\relax\else
\ifnum#1>\maxEqNo\global\maxEqNo=#1\fi%
\expandafter\ifx\csname E#1\endcsname\FTEST\relax\else
\immediate\write16{E#1 hat sich geaendert!}\fi
\expandwrite\AUX{\neverexpand\ref{E#1}{\FTEST}}\fi
\llap{\hbox to 40pt{\marginnote{#1}\FTEST\hfill}}}

\catcode`@=11
\def\eqalignno#1{\null\vcenter{\openup\jot\m@th\ialign{\eqno{##}\hfil
&\strut\hfil$\displaystyle{##}$&$\displaystyle{{}##}$\hfil\crcr#1\crcr}}\,}
\catcode`@=12

%%%%%%%%%%%%
\newbox\QEDbox
\newbox\nichts\setbox\nichts=\vbox{}\wd\nichts=2mm\ht\nichts=2mm
\setbox\QEDbox=\hbox{\vrule\vbox{\hrule\copy\nichts\hrule}\vrule}
\def\qed{\leavevmode\unskip\hfil\null\nobreak\hfill\copy\QEDbox\medbreak}
%%%%%%%%%%%%%%
\newdimen\HIindent
\newbox\HIbox
\def\setHI#1{\setbox\HIbox=\hbox{#1}\HIindent=\wd\HIbox}
\def\HI#1{\par\hangindent\HIindent\hangafter=0\noindent\leavevmode
\llap{\hbox to\HIindent{#1\hfil}}\ignorespaces}
%%%%%%%%%%%%%%

\newdimen\maxSpalbr
\newdimen\altSpalbr
\newcount\Zaehler

%{\catcode`/=\active
%\gdef\SlashOn{\catcode`/=\active\def/{X\string/\hskip0pt Y}}
%}

\newif\ifxxx

{\catcode`/=\active

\gdef\beginrefs{%
\xxxfalse
\catcode`/=\active
\def/{\string/\ifxxx\hskip0pt\fi}
\def\TText##1{{\xxxtrue\tt##1}}
\expandafter\ifx\csname Spaltenbreite\endcsname\relax
\def\Spaltenbreite{1cm}\immediate\write16{Spaltenbreite undefiniert!}\fi
\expandafter\altSpalbr\Spaltenbreite
\maxSpalbr0pt
\gdef\alt{}
\def\\##1\relax{%
\gdef\neu{##1}\ifx\alt\neu\global\advance\Zaehler1\else
\xdef\alt{\neu}\global\Zaehler=1\fi\xdef\SigText{##1\the\Zaehler}}
\def\L|Abk:##1|Sig:##2|Au:##3|Tit:##4|Zs:##5|Bd:##6|S:##7|J:##8|xxx:##9||{%
\def\SigText{##2}\global\setbox0=\hbox{##2\relax}
\edef\TEST{[\SigText]}
\expandafter\ifx\csname##1\endcsname\TEST\relax\else
\immediate\write16{##1 hat sich geaendert!}\fi
\expandwrite\AUX{\neverexpand\ref{##1}{\TEST}}
\setHI{[\SigText]\ }
\ifnum\HIindent>\maxSpalbr\maxSpalbr\HIindent\fi
\ifnum\HIindent<\altSpalbr\HIindent\altSpalbr\fi
\HI{\marginnote{##1}[\SigText]}
\ifx-##3\relax\else{##3}: \fi
\ifx-##4\relax\else{##4}{\sfcode`.=3000.} \fi
\ifx-##5\relax\else{\it ##5\/} \fi
\ifx-##6\relax\else{\bf ##6} \fi
\ifx-##8\relax\else({##8})\fi
\ifx-##7\relax\else, {##7}\fi
\ifx-##9\relax\else, \TText{##9}\fi\Par}
\def\B|Abk:##1|Sig:##2|Au:##3|Tit:##4|Reihe:##5|Verlag:##6|Ort:##7|J:##8|xxx:##9||{%
\def\SigText{##2}\global\setbox0=\hbox{##2\relax}
\edef\TEST{[\SigText]}
\expandafter\ifx\csname##1\endcsname\TEST\relax\else
\immediate\write16{##1 hat sich geaendert!}\fi
\expandwrite\AUX{\neverexpand\ref{##1}{\TEST}}
\setHI{[\SigText]\ }
\ifnum\HIindent>\maxSpalbr\maxSpalbr\HIindent\fi
\ifnum\HIindent<\altSpalbr\HIindent\altSpalbr\fi
\HI{\marginnote{##1}[\SigText]}
\ifx-##3\relax\else{##3}: \fi
\ifx-##4\relax\else{##4}{\sfcode`.=3000.} \fi
\ifx-##5\relax\else{(##5)} \fi
\ifx-##7\relax\else{##7:} \fi
\ifx-##6\relax\else{##6}\fi
\ifx-##8\relax\else{ ##8}\fi
\ifx-##9\relax\else, \TText{##9}\fi\Par}
\def\Pr|Abk:##1|Sig:##2|Au:##3|Artikel:##4|Titel:##5|Hgr:##6|Reihe:{%
\def\SigText{##2}\global\setbox0=\hbox{##2\relax}
\edef\TEST{[\SigText]}
\expandafter\ifx\csname##1\endcsname\TEST\relax\else
\immediate\write16{##1 hat sich geaendert!}\fi
\expandwrite\AUX{\neverexpand\ref{##1}{\TEST}}
\setHI{[\SigText]\ }
\ifnum\HIindent>\maxSpalbr\maxSpalbr\HIindent\fi
\ifnum\HIindent<\altSpalbr\HIindent\altSpalbr\fi
\HI{\marginnote{##1}[\SigText]}
\ifx-##3\relax\else{##3}: \fi
\ifx-##4\relax\else{##4}{\sfcode`.=3000.} \fi
\ifx-##5\relax\else{In: \it ##5}. \fi
\ifx-##6\relax\else{(##6)} \fi\PrII}
\def\PrII##1|Bd:##2|Verlag:##3|Ort:##4|S:##5|J:##6|xxx:##7||{%
\ifx-##1\relax\else{##1} \fi
\ifx-##2\relax\else{\bf ##2}, \fi
\ifx-##4\relax\else{##4:} \fi
\ifx-##3\relax\else{##3} \fi
\ifx-##6\relax\else{##6}\fi
\ifx-##5\relax\else{, ##5}\fi
\ifx-##7\relax\else, \TText{##7}\fi\Par}
\bgroup
\baselineskip12pt
\parskip2.5pt plus 1pt
\hyphenation{Hei-del-berg Sprin-ger}
\sfcode`.=1000
\beginsection References. References

}}

\def\endrefs{%
\expandwrite\AUX{\neverexpand\ref{Spaltenbreite}{\the\maxSpalbr}}
\ifnum\maxSpalbr=\altSpalbr\relax\else
\immediate\write16{Spaltenbreite hat sich geaendert!}\fi
\egroup\write16{Letzte Gleichung: E\the\maxEqNo}}

%\L|Abk:|Sig:|Au:|Tit:|Zs:|Bd:|S:|J:|xxx:-||
%\B|Abk:|Sig:|Au:|Tit:|Reihe:|Verlag:|Ort:|J:|xxx:-||
%\Pr|Abk:|Sig:|Au:|Artikel:|Titel:|Hgr:|Reihe:|Bd:|Verlag:|Ort:|S:|J:|xxx:-||

\let\hut\^
\def\hochoplus#1{^{\oplus#1}}
\def\^{\ifmmode\expandafter\hochoplus\else\expandafter\hut\fi}

\def\notdiv{\mathrel{\not\kern-.5pt|}}
\def\otimes{\Otimes}

\catcode`@=11
\def\Ddots{\mathinner{
\mkern1mu\raise0\p@\vbox{\kern0\p@\hbox{.}}
\mkern2mu\raise4\p@\hbox{.}
\mkern2mu\raise8\p@\hbox{.}
\mkern1mu}}
\def\DDdots{\mathinner{
\mkern1mu\raise1.45\p@\vbox{\kern0\p@\hbox{.}}
\mkern2mu\raise4\p@\hbox{.}
\mkern2mu\raise6.55\p@\hbox{.}
\mkern1mu}}
\catcode`@=12

%\marginalnotestrue
\marginalwarningsfalse

\title{On Noether's and Weyl's bound in positive characteristic}

{\narrower

\noindent{\bf Abstract:} In this note we generalize several well known
results concerning invariants of finite groups from characteristic
zero to positive characteristic not dividing the group order. The
first is Schmid's relative version of Noether's theorem. That theorem
compares the degrees of generators of a group with those of a
subgroup. Then we prove a suitable positive characteristic version of
Weyl's theorem on vector invariants: polarization works in small
degrees. Using that we show that the regular representation has the
``most general'' ring of invariants, thereby generalizing theorems of
Schmid and Smith.

}

\beginsection Intro. Introduction

Let $G$ be a finite group acting linearly on a $k$\_vector space
$V$. If $k=\CC$, Noether showed in \cite{Noe} that the ring $k[V]^G$
of $G$\_invariants is generated by invariants of degree $\le|G|$. This
is called {\it Noether's bound}. An inspection of her proof reveals
that it stays valid for any field $k$ of characteristic~0 or
characteristic $p>|G|$.  It has been a long standing question whether
Noether's bound holds under the weaker assumption that $p$ does not
divide $|G|$. This has been proved only recently by Fleischmann
\cite{Fl} and Fogarty \cite{Fo}.

Meanwhile, Schmid showed in \cite{Sch} that Noether's bound is
sharp only for cyclic groups. Again, her proof is valid if either
$p=0$ or $p>|G|$. The only reason why Schmid's proof doesn't carry
over to characteristics not dividing $|G|$ is a certain refinement of
Noether's theorem which compares degrees of generators of
$G$\_invariants with those of $H$\_invariants where $H$ is any
subgroup of $G$. The purpose of this paper is to provide this tool for
$p\notdiv|G|$ (see \cite{TD} for a precise statement).

After a first draft of this note has been completed it came to my
attention that this has also been done independently by Fleischmann,
\cite{Fl}, and Sezer, \cite{Se}. Nevertheless, we find it still
worthwhile to include a proof. First, our proof is an adaptation of Fogarty's
(as opposed to Fleischmann's) proof of Noether's bound. Secondly, we
introduce another numerical invariant of independent interest. More
precisely, in section~\cite{NullCone} we look at the zero fiber of the
quotient map $V\pfeil V/G$ (the null cone) which is the spectrum
of a local Artinian ring. We study the degree of nilpotence of its
maximal ideal in \cite{TA} and show that it is closely related to
degrees of generating invariants.

In a further section, we study degree bounds for algebras over
arbitrary ground rings $k$. We show that it suffices to consider for
$k$ a field and even there the problem depends only on the
characteristic of $k$ with $\|char|k=0$ giving the lowest bounds.

The second major thread of this paper is Weyl's theorem on vector
invariants. It asserts that one can obtain all invariants on any
number of copies of a representation $V$ from invariants which live on
only $\|dim|V$ copies. The process allowing this is called
polarization and works in that generality only in characteristic
zero. Nevertheless, in positive characteristic it is plausible that
polarization still works in ``low'' degrees. \cite{Polar} makes this
statement precise. On the other hand, generating invariants for finite
groups have ``low'' degree by the theorem of
Noether\_Fleischmann\_Fogarty. This yields Weyl type theorems for
vector invariants in positive characteristic. In particular, we
strengthen theorems of Schmid and Smith according to which the regular
representation has the most general ring of invariants.

It should be added that even in characteristic zero our approach is
not without interest. It yields a completely elementary proof
of Weyl's theorem for the (special but most important) case that $G$
is linearly reductive. In particular, we avoid the use of Capelli's
identity and representation theory of $GL(n)$.

\medskip\noindent {\bf Acknowledgment:} The author would like to thank
the referee for very carefully reading the manuscript.

\beginsection NullCone. The degree of nilpotence of the null cone

Let $A$ be a commutative ring. For a subset $S$ let $\<S\>$ be the
ideal generated by $S$. Now assume that the finite group $G$ acts on
$A$. For a $G$\_invariant ideal $\fm\subseteq A$ we define
$$\eqno{}
\eta_G(\fm):=\|inf|\{d\in\NN\mid\fm^d\subseteq\<\fm^G\>\}.
$$
It is well known that for finite groups, the fibers of
$\|Spec|A\pfeil\|Spec|A^G$ are precisely the $G$\_orbits. Thus $\fm$
and $\<\fm^G\>$ define the same vanishing set in $\|Spec|A$, which
implies $\eta_G(\fm)<\infty$ whenever $A$ is Noetherian (see also
\cite{CZ} below). The most important case will be when $A$ is the ring
of polynomial functions on a representation of $G$ and $\fm$ the
maximal ideal corresponding to the origin. Then $\<\fm^G\>$ defines
the null cone and $\eta_G(\fm)-1$ is the degree of nilpotence of the
maximal ideal of $A/\<\fm^G\>$.

\Theorem TA. Let $H\subseteq G$ be a subgroup and assume $[G:H]\in
A^\Times$. Then
$$\eqno{}
\eta_G(\fm)\le[G:H]\,\eta_H(\fm).
$$

\Proof: Let $d:=[G:H]$. Choose $d$ arbitrary elements $a_u\in \fm^H$
which are indexed by $u\in G/H$. In the following, we will regard $u$ as an element of $G$ representing its coset. Consider the expression
$$\eqno{}
\Phi:=\sum_{u\in G/H}\prod_{v\in G/H}(va_v-ua_v)
$$
Since every product contains a zero factor, namely the one corresponding to
$v=u$, we have $\Phi=0$. On the other side, we can expand the product
and get
$$\eqno{}
\Phi=\sum_S(-1)^{|S|}\Phi_S
$$
where $S$ ranges over all subsets of $G/H$ and
$$\eqno{}
\Phi_S=\sum_{u\in G/H}\prod_{v\not\in S}(va_v)\prod_{v\in S}(ua_v)=
\prod_{v\not\in S}(va_v)
\left[\sum_{u\in G/H}u\Big(\prod_{v\in S}a_v\Big)\right].
$$
Thus, $\Phi_S\in\<\fm^G\>$ unless $S=\emptyset$ when
$\Phi_\emptyset=d\,\prod_{v\in G/H}va_v$. Since the $a_v$ were
arbitrary and $d\in A^\Times$ we get
$$\eqno{}
\prod_{v\in G/H}v\<\fm^H\>\subseteq\<\fm^G\>.
$$
With $e:=\eta_H(\fm)$ we have by definition
$\fm^e\subseteq\<\fm^H\>$. Thus,
$$\eqno{}
\<\fm^G\>\supseteq\prod_{v\in G/H}v\fm^e=\fm^{de}.
$$
This implies $\eta_G(\fm)\le de=[G:H]\,\eta_H(\fm)$.\qed

\noindent For $H=1$ we get:

\Corollary CZ. Assume $|G|\in A^\Times$. Then $\eta_G(\fm)\le|G|$.

\beginsection degrees. The degree of generators

Let $A=\oplus_{d=0}^\infty A_d$ be an $\NN$\_graded commutative
ring. We use the notation $A_{\le n}:=\oplus_{d=0}^n A_d$ and
analogously $A_{\ge n}:=\oplus_{d=n}^\infty A_d$. Of particular
interest is the ideal $A_{\ge1}$ which we denote by $\fm$. We define
$$\eqno{}
\beta(A):=\|inf|\{n\in\NN\mid A\hbox{ is generated as a ring by }A_{\le n}\},
$$
which is the smallest degree of a generator of $A$. Now assume that the finite
group $G$ acts on $A$ by degree preserving automorphisms.

\Lemma TB. Assume $\beta(A)\le1$ and $|G|\in A^\Times$. Then
$\beta(A^G)\le\eta_G(\fm)$.

\Proof: Let $e:=\eta_G(\fm)$. By the classical argument of Hilbert,
any set of $G$\_invariant generators of the ideal $\<\fm^G\>$ also
generates the ring $A^G$. Thus, we have to show that $\<\fm^G\>$ is
generated by invariants of degree $\le e$. To that end, it suffices to
show that $\<\fm^G\>$ is generated by elements (invariant
or not) of degree $\le e$ since each of them can be expressed in terms
of invariants of degree $\le e$.

By definition, $\fm^e\subseteq\<\fm^G\>$. Since $\beta(A)\le 1$ we
have $\fm^e=A_{\ge e}$ and that ideal is generated by $A_e$. Thus,
$\<\fm^G\>$ is generated by $A_{\le e-1}^G$ and $A_e$. This implies
that $\<\fm^G\>$ is generated by $A_{\le e}^G$.\qed

For the next statement let $A(G)$ be the polynomial ring over $A$ in
$|G|$ variables $\{x_g\mid g\in G\}$ with the $G$\_action extended by
$gx_h:=x_{gh}$. It is also $\NN$\_graded with $\|deg|x_g=1$. We can
write $A(G)=A\otimes_\ZZ\ZZ(G)$ where $\ZZ(G)=S^*(\ZZ[G])$ is the
symmetric algebra over the (integral) regular representation of
$G$. Using this, we see that $A(G)$ is actually bigraded where $a\in
A_d$ and $x_g$ have bidegree $(d,0)$ and $(0,1)$, respectively. The
following result is a partial converse of \cite{TB}.

\Lemma TC. For any $A$ we have that $\eta_G(\fm)\le\beta(A(G)^G)$.

\Proof: Let $b:=\beta(A(G)^G)$. For $a\in A_e$
with $e\ge b$ let
$$\eqno{}
\aq:=\sum_{u\in G}u(ax_1)=\sum_{u\in G}(ua)x_u\in A(G)^G.
$$
Its bidegree is $(e,1)$. Since $A(G)^G$ is generated by elements of
degree $\le b$ we can write $\aq$ as a sum of products $p\cdot q$
where $p$ and $q$ are bihomogeneous $G$\_invariants of bidegree
$(c,0)$ and $(e-c,1)$ respectively. We may assume that
$\|deg|q=e-c+1\le b$. This implies $c\ge e-b+1\ge1$. Thus,
$p\in\fm^G$. Comparing the coefficient of $x_1$ in the expression
$\aq=\sum pq$, we obtain $a\in \<\fm^G\>$. Thus $\fm^b\subseteq A_{\ge
b}\subseteq\<\fm^G\>$.\qed

\noindent Combining these comparison lemmas with \cite{TA}, we obtain:

\Theorem TD. Assume $\beta(A)\le1$ and $|G|\in A^\Times$. Let $H\subseteq
G$ be a subgroup. Then
$$\eqno{}
\beta(A^G)\le[G:H]\beta(A(H)^H).
$$

\Proof: Indeed,
$$\eqno{}
\vcenter{\halign{$#\,$\hfill&$\,#$\hfill&\qquad (by \cite{#})\hfill\cr
\beta(A^G)&\le\eta_G(\fm)&TB\cr
&\le[G:H]\eta_H(\fm)&TA\cr
&\le[G:H]\beta(A(H)^H)&TC\cr}}
$$\qed

\Remark: Since \cite{TA} holds under the weaker assumption
$[G:H]\in A^\Times$ one might wonder whether also \cite{TD} holds in
that generality.
\medskip

\noindent
Taking for $H$ the trivial subgroup, we obtain Fleischmann's and
Fogarty's theorem:

\Corollary CA. Assume $\beta(A)\le1$ and $|G|\in A^\Times$. Then
$\beta(A^G)\le|G|$.

\noindent This suggests a definition. For a fixed ground ring $k$
let $\beta_k(G):=\sup_A\beta(A^G)$ where $A$ runs through all
$\NN$\_graded $k$\_algebras with $G$\_action and with $\beta(A)\le1$.
The following generalizes a theorem of Schmid~\cite{Sch} for $k=\CC$:

\Corollary CB. Assume $|G|\in k^\Times$ and let $H\subseteq G$ be a
subgroup. Then
$$\eqno{}
\beta_k(G)\le[G:H]\beta_k(H).
$$
In particular, $\beta_k(G)\le|G|$.

\def\AS{{\tilde A}}
For arbitrary $\NN$\_graded algebras we have

\Theorem TE. Assume $|G|\in k^\Times$ and let $A$ be a $k$\_algebra. Then
$\beta(A^G)\le\beta_k(G)\beta(A)$.

\Proof: Let $b:=\beta(A)$ and let $\Aq$ be the free commutative
$A_0$\_algebra over the $A_0$\_module $\oplus_{d=1}^bA_d$. Then $\Aq$
has two degree functions, $\|deg|$ and $\|deg|'$, where an element of
$A_d$ has degree $d$ or $1$, respectively. Clearly, we have $\|deg|\le
b\|deg|'$. If we denote $\Aq$ equipped with the grading $\|deg|'$ by
$\AS$ then $\beta(\AS)\le1$.

There is a natural map of $\Aq$ onto $A$. Hence $\Aq^G\pfeil A^G$ is
also onto. Therefore,
$$\eqno{}
\beta(A^G)\le\beta(\Aq^G)\le b\,\beta(\AS^G)\le b\, \beta_k(G).
$$\qed

The following consequence is due to Schmid \cite{Sch} for $k=\CC$ and
to Fleischmann \cite{Fl} in general:

\Corollary CD. Let $N\triangleleft G$ be a normal
subgroup with $[G:N]\in k^\Times$ and $A$ a $k$\_algebra. Then
$$\eqno{}
\beta(A^G)\le\beta_k(G/N)\beta(A^N)
$$
In particular,
$\beta_k(G)\le\beta_k(G/N)\beta_k(N)$.

\Proof: Apply \cite{TE} to $G/N$ acting on $A^N$.\qed

Now, we can state the following generalization of Schmid's result:

\Theorem TF. Let $G$ be a non-cyclic group with $|G|\in k^\Times$. Then
$\beta_k(G)<|G|$.

\Proof: A close inspection shows that Schmid's proof works now (i.e.,
after having proved \cite{CB}) for any $\NN$\_graded $k$\_algebra.\qed

For completeness we state two more results without proof. In
characteristic zero, Domokos and Heged\H us strengthened Schmid's
theorem to:

\Theorem DomHeg. {\rm(\cite{DH})} Assume $k=\CC$. Let
$G$ be a non\_cyclic finite group. Then
$$\eqno{}
\beta_k(G)\le\cases{{3\over 4}|G|&if $|G|$ is even\cr
{5\over8}|G|&if $|G|$ is odd\vrule height 12pt
width 0pt\relax\cr}
$$

\noindent The only part of the proof which cannot be easily carried
over to positive characteristic is Schmid's calculation of $\beta_k$
for dihedral groups. This has been done by Sezer who thereby
proved

\Theorem DomHegSez. {\rm(\cite{Se})} \cite{DomHeg} holds for any field
$k$ whose characteristic does not divide $|G|$.

In analogy to $\beta_k(G)$ we can define $\eta_k(G)$ as
$\|sup|_{A,\fm}\eta_G(\fm)$ where $A$ runs through all (not
necessarily graded) $k$\_algebras with $G$\_action and $\fm$ through
all $G$\_invariant ideals of $A$. Then \cite{CZ} says
$\eta_k(G)\le|G|$ whenever $|G|\in k^\Times$. In fact, we have:

\Theorem etabeta. For all finite groups $G$ with $|G|\in k^\Times$ holds
$\eta_k(G)=\beta_k(G)$.

\Proof: \cite{TB} shows $\beta_k(G)\le\eta_k(G)$. For the
converse inequality, let $\fm$ be a $G$\_invariant ideal in the
$G$\_ring $A$. Let $\AS$ be the blow-up algebra of $A$ along $\fm$,
i.e., $\AS$ is the $\NN$\_graded subring $\oplus_{d=0}^\infty\fm^dt^d$
of $A[t]$. Let $\tilde\fm:=\AS_{\ge1}$. Then \cite{TC} implies
$$\eqno{}
\eta_G(\tilde\fm)\le\beta(\AS(G)^G)\le\beta_k(G).
$$
Let $\pi:\AS\pfeil A$ be the composition of $\AS\into A[t]\auf A$
where $t$ maps to $1$. By linear reductivity, we have
$\pi(\<\tilde\fm^G\>)=\<\fm^G\>$ and (obviously)
$\pi(\tilde\fm^d)=\fm^d$. This implies
$\eta_G(\fm)\le\eta_G(\tilde\fm)$, hence
$\eta_k(G)\le\beta_k(G)$.\qed

\Remark: If $\|char|k$ divides $|G|$ then Richman has shown
$\beta_k(G)=\infty$ (see \cite{modular} below). It is feasible that
this is also true for $\eta_k(G)$. Then one could remove the
assumption $|G|\in k^\Times$ in \cite{etabeta}. As an example, we
compute $\eta_k(G)$ where $G=\{1,\sigma\}$ is the group of order $2$
and $\|char|k=2$. Let $V=k^{2n}$ with coordinates
$x_1,\ldots,x_n,y_1,\ldots,y_n$ and assume $\sigma(x_i)=y_i$. Let
$A=k[V]$. Then $A^G$ is generated by all $x_iy_i$, $1\le i\le n$ and
$x^\alpha+y^\alpha$, $\alpha\in\NN^n$. In particular, $\fm^G$ contains
$x_i+y_i$. That means $A/\<\fm^G\>\cong k[x_1,\ldots,x_n]/I$ where $I$
is the ideal generated by all $-x_i^2$ and all
$(1+(-1)^{|\alpha|})x^\alpha$ with $|\alpha|>0$. In characteristic
two, the second kind of elements are all zero. Therefore,
$A/\<\fm^G\>=k[x_1,\ldots,x_n]/(x_1^2,\ldots,x_n^2)$. Since $x_1\ldots
x_n\not\in I$ we have $\eta_G(\fm)=n+1$. Thus, $\eta_k(G)=\infty$.

\beginsection Functoriality. Some functoriality properties

In this section, we develop some machinery to compute $\beta_k(A)$ for
arbitrary ground rings $k$.

\Theorem faith. Let $\phi:k\pfeil K$ be a ring homomorphism. Then
$\beta_k(G)\ge\beta_K(G)$ with equality if $\phi$ is faithfully
flat.

\Proof: Let $A$ be a graded $K$\_algebra with $G$\_action and
$\beta(A)\le1$. By means of $\phi$, we can consider $A$ as a
$k$\_algebra. Since $\beta(A)$ and $\beta(A^G)$ depend only on the
underlying ring structure, we get $\beta_k(G)\ge\beta_K(G)$.

Now assume that $\phi$ is faithfully flat and let $A$ be a graded
$k$\_algebra with $G$\_action and $\beta(A)\le1$. Put $B:=A^G$,
$\Aq:=A\otimes_kK$ and $\Bq:=B\otimes_kK$. Then $\Bq=\Aq^G$. Indeed,
by flatness, the exact sequence
$$\eqno{}
0\Pfeil B\Pfeil A\Pf{ga-a}\oplus_{g\in G}A
$$
stays exact upon tensoring with $K$. Now we claim
$\beta(\Bq)\ge\beta(B)$. For $b:=\beta(\Bq)$ and any $d\in\NN$
let $C$ be the cokernel of
$$\eqno{5}
\bigoplus_{\sum i_j=d\atop 0\le i_j\le b}B_{i_1}\otimes\ldots\otimes B_{i_s}
\Pfeil B_d.
$$
This map becomes surjective after tensoring with $K$. Thus,
$C\otimes_kK=0$. Faithful flatness implies $C=0$ and the claim follows.

We conclude $\beta(A^G)\le\beta(\Aq^G)\le\beta_K(G)$ for arbitrary
$A$, hence $\beta_k(G)\le\beta_K(G)$.\qed

\Corollary. Let $k$ be a field. Then $\beta_k(A)$ depends only on the
characteristic of $k$.

\noindent In view of this result, we write $\beta_p(G):=\beta_k(G)$
where $k$ is any field of characteristic $p$.

Next, we cite a result of Richman,~\cite{Ri}, which shows that in the
modular case there is no universal bound on the degree of generators:

\Theorem modular. Let $k$ be a field whose characteristic divides
$|G|$. Then $\beta(S^*_k(U)^G)$ is unbounded, where $U$ runs through
all finitely generated $k[G]$\_modules. In particular,
$\beta_k(G)=\infty$.

\noindent This is being used in the proof of:

\Proposition vector. Let $k$ be an arbitrary ring. Then, to compute
$\beta_k(G)$ it suffices to take the supremum over all
$\beta(S^*_k(U)^G)$ where $U$ is a finitely generated $k[G]$\_module.

\Proof: First assume that $|G|$ is not invertible in $k$. Then there
is a prime divisor $p$ of $|G|$ with $p\not\in k^\Times$, i.e., $pk\ne
k$. Let $\fm\subset k$ be a maximal ideal containing $pk$ and
$K:=k/\fm$. Then $\|char|K=p$ and the assertion follows from
\cite{modular}.

Thus, we may assume $|G|\in k^\Times$. Let $A$ be any graded
$k$\_algebra with $G$\_action and $\beta(A)=1$. Let
$b:=\beta(A^G)$. Then there is $f\in A_b^G$ which is not a polynomial
of lower degree invariants. Clearly, we may assume $b\ge1$. Hence, $f$
is a polynomial of finitely many elements of $A_1$ with coefficients
in $k$. Let $U$ be the $k[G]$\_module generated by these
elements. Then we get an algebra homomorphism $S^*(U)\pfeil A$ whose
image contains $f$. Since $G$ is linearly reductive, we can lift $f$
to an invariant $\tilde f\in S^b_k(U)^G$. Clearly, we can't obtain
$\tilde f$ as a polynomial of lower degree invariants, either. Hence
$\beta(S^*(U)^G)\ge\beta(A^G)$ which shows the assertion.  \qed

Now we can reduce the computation of $\beta_k(G)$ for an arbitrary
ring $k$ to fields. For this we define
$$\eqno{}
\|char|(k):=\cases{\{0\}&if $\QQ\subseteq k$;\cr
\{p\mid p\not\in k^\Times\}&otherwise.\cr}
$$

\Theorem generalk. Let $k$ be any ring. Then
$$\eqno{6}
\beta_k(G)=\|sup|\{\beta_p(G)\mid p\in\|char|(k)\}.
$$

\Proof: If $\QQ\subseteq k$ then $k$ is faithfully flat over $\QQ$ and
the assertion follows from \cite{faith}. Thus, assume
$\QQ\not\subseteq k$. If $p$ is not invertible in $k$ then there is a
maximal ideal $\fm$ such that $k/\fm$ has characteristic $p$. Thus,
$\beta_k(G)$ is at least as big as the left hand side of \cite{E6}. In
particular, we are done if one of the $\beta_p(G)$ is $\infty$. Thus
we may assume that $|G|$ is invertible in $k$.

Let $A$ be a $k[G]$\_algebra which, by \cite{vector}, we may assume to
be of the form $A=S^*_k(U)$ where $U$ is a finitely generated
$k[G]$\_module. This implies that each homogeneous component $A_d$ is
a finitely generated $k$\_module. Being a direct summand the
same holds for $A^G_d$.

Let $\fm\subset k$ be a maximal ideal. Then we claim that
$$\eqno{8}
A^G\otimes_kk/\fm\pfeil(A\otimes_kk/\fm)^G
$$
is an isomorphism. Indeed, injectivity follows from the fact that $A^G$ is a
direct summand of $A$. Moreover,
$A^G\pfeil(A/\fm A)^G=(A\otimes_kk/\fm)^G$ is surjective since $|G|$ is
invertible in $k$.

Let $b$ be the right hand side of \cite{E6} and let
$C(b,d)$ be the cokernel of \cite{E5}. Then $C(b,d)\otimes_kk/\fm=0$
for every maximal ideal. Since $C(b,d)$ is finitely generated this
implies $C(b,d)=0$ (Nakayama Lemma). Thus $\beta(A^G)\le b$ and
therefore $\beta_k(G)\le b$.\qed

\Corollary. For any ring $k$ holds $\beta_k(G)<\infty$ if and only if
$|G|\in k^\Times$.

A similar argument shows that $k=\QQ$ is the best possible case:

\Theorem. For all primes $p$ holds $\beta_p(G)\ge\beta_0(G)$. For
almost all primes holds equality.

\Proof: If $p$ divides $|G|$ there is nothing to prove since then
$\beta_p(G)=\infty$ (\cite{modular}). Thus assume $p\notdiv|G|$. Let
$k=\ZZ_p$ be the ring of $p$\_adic integers with field of fractions
$K=\QQ_p$ and residue field $\FF_p$. By \cite{vector}, there is a
finite dimensional $K[G]$\_module $U_K$ with
$\beta(S^*_K(U_K)^G)=\beta_0(G)=:b_0$. Let $U_k\subseteq U_K$ be a
$G$\_stable lattice and $A:=S^*_k(U_k)$. Let
$\Aq:=A\otimes_k\FF_p$. Since $\Aq$ is an $\FF_p$\_algebra, it
suffices to prove $\beta(\Aq^G)\ge b_0$.

The same argument as for \cite{E8} shows that
$$\eqno{7}
A^G\otimes_k\FF_p\pfeil(A\otimes_k\FF_p)^G
$$
is an isomorphism. Now let $B:=A^G$ and consider for any $b$ again the
cokernel $C(b,d)$ of the map \cite{E5}. For any $b<b_0$ there is a
$d>0$ such that $C(b,d)\otimes_kK\ne0$. But $C(b,d)$ is a finitely
generated $k$\_module which implies $C(b,d)\otimes_k\FF_p\ne0$. This
implies that $\beta(\Aq^G)$ can not be smaller than $b$.

Now we show that conversely $\beta_p(G)\le\beta_0(G)$ for almost all
primes. Let $k(G)$ be the symmetric algebra over the regular
$k$\_representation. Let $f_1,\ldots,f_r$ be a minimal set of
generators of $\QQ(G)^G$. Let $h_1,\ldots,h_s$ be generators of
$\ZZ[{1\over g}](G)^G$ where $g:=|G|$. Then there is a multiple $N$ of
$g$ such that the $f_i$ are defined over $k:=\ZZ[{1\over N}]$ and such
that all $h_j$ are in the $k$\_algebra generated by the $f_i$. This
implies that $k(G)^G$ is also generated by the $f_i$. The homomorphism
$k(G)^G\pfeil\FF_p(G)^G$ is surjective where $p$ is any prime not
dividing $N$. This implies $\beta(\FF_p(G)^G)\le\beta(\QQ(G)^G)$.  We
conclude with theorems of Schmid~\cite{Sch} and Smith~\cite{Sm} (see
also its generalization, \cite{regrep2}, further down) which assert
$\beta(\QQ(G)^G)=\beta_0(G)$ and $\beta(\FF_p(G)^G)=\beta_p(G)$ for
$p>|G|$.\qed

\Remark: Presently, no group $G$ and prime $p$ not
dividing $|G|$ with $\beta_p(G)>\beta_0(G)$ seems to be known.

\beginsection Polarization. Polarization in positive characteristic

From now on, $k$ will be a field which, for convenience, we assume to
be infinite. We set $p=\|char|k$ if positive and $p=\infty$
otherwise. Let $V$ be a $k$\_vector space of dimension
$\ell<\infty$. Let $k[V]_d=S^d(V^*)$ be the space of polynomial
functions which are homogeneous of degree $d$. We are interested in
the space $V\^n$ of $n$\_tuples of vectors in $V$.

There is a natural action of the algebraic group $GL(n)$ on
$V\^n=V\otimes_kk^n$. Let $S$ be a subset of $k[V\^n]$. Then the
$GL_n$\_submodule $\<S\>_{GL(n)}$ generated by $S$ is called the {\it
polarization\/} of $S$. The action of the Lie algebra is given by the
{\it polarization operators\/}
$$\eqno{}
P_{jj'}:=\sum_{i=1}^\ell x_{ij}{\partial\over\partial x_{ij'}},
\qquad 1\le j,j'\le n
$$
where $x_{1j},\ldots,x_{\ell j}$ are coordinates of the $j$-th copy of
$V$ in $V\^n$. If $m\le n$ then we consider $V\^m$ as
a quotient of $V\^n$ by forgetting the last $n-m$ components.
In characteristic zero, it follows from work by Weyl (\cite{Weyl}~II.5)
that every function on $V^n$ with $n\ge\ell:=\|dim|V$ can be obtained
by polarization from functions on $V^\ell$. In positive characteristic
this is true in sufficiently small degrees.

\Theorem Polar. Let $n\ge m\ge\ell:=\|dim|V$ and $d\le(p-1)m$. Then
$$\eqno{}
k[V\^n]_d=\<k[V\^m]_d\>_{GL(n)}.
$$

\Proof: First, by replacing $V$ with $V\oplus k^{m-\ell}$ we may
assume that $m=\ell$. Moreover, it clearly suffices to treat the case
$n=m+1=\ell+1$. The proof will proceed by induction on $\ell$ starting
at the trivial case $\ell=0$.

We label the copies of $V\^{\ell+1}$ by $j=0,\ldots,\ell$ such that
$V\^\ell$ corresponds to $j=1,\ldots,\ell$. The coordinates of the
$j$\_th copy are denoted by $x_{1j},\ldots,x_{\ell j}$. A monomial in
the $x_{ij}$ can be represented by an $\NN$\_valued matrix
$$\eqno{1}
A=\pmatrix{a_{10}&a_{11}&\ldots&a_{1\ell}\cr
\vdots&\vdots&\Ddots&\vdots\cr
a_{\ell0}&a_{\ell1}&\ldots&a_{\ell\ell}\cr}
$$
where $x^A:=\prod_{ij}x_{ij}^{a_{ij}}$. Its degree is
$\|deg|A:=\sum_{ij}a_{ij}$.

We introduce a total order on the set of matrices $A$ as
follows. Let $r_i=r_i(A)=\sum_j{a_{ij}}$ and $c_j=c_j(A)=\sum_ia_{ij}$
be the row and column sums, respectively. Then we define the {\it
index} of $A$ as the vector
$$\eqno{}
\|ind|A=(c_0,\ldots,c_\ell,r_1,\ldots,r_\ell,a_{10},\ldots,a_{1\ell},
a_{20},\ldots,a_{2\ell},\ldots,a_{\ell0},\ldots,a_{\ell\ell}).
$$
In words: $\|ind|A$ starts with the column sums
followed by the row sums followed by the coefficients read from left
to right, top to bottom. We say $A < B$ if $\|ind|A$ is
lexicographically smaller than $\|ind|B$. We have to show that
$\|deg|A=d\le(p-1)\ell$ implies $x^A\in
X(\ell,d):=\<k[V\^\ell]_d\>_{GL(\ell+1)}$. Let $A$ be the
smallest counterexample.

First, observe that $X(\ell,d)$ is invariant under both $GL(\ell)$
(acting on $V$) and $GL(\ell+1)$ (by definition). Thus, the property
$x^A\in X(\ell,d)$ is invariant under permutations of the rows or the
columns of $A$. The minimality of $A$ implies
$$
\eqalignno{&&c_0(A)\le\ldots\le c_\ell(A)\quad\hbox{and}\cr
&&r_1(A)\le\ldots\le r_\ell(A).\cr}
$$
From $d=\sum_ir_i(A)\le(p-1)\ell$ we obtain
$$\eqno{2}
r_1(A)\le p-1.
$$

The matrix $A$ must have a non\_zero entry in the first column since
otherwise $x^A\in k[V\^\ell]_d\subseteq X(\ell,d)$. A
fortiori, there is non\_zero entry $a_{ij}$ with $i+j\le\ell$ (i.e.,
which is strictly above the dotted diagonal in \cite{E1}). Let
$a_{ij}$ be the first one, where we read from left to right and top to
bottom. Then $a_{ij}$ is not in the last column and the entries above
$a_{ij{+}1}$ are all zero.

Let $\Aq$ be the matrix obtained from $A$ by replacing the
entries $a_{ij}$ and $a_{ij{+}1}$ by $a_{ij}-1$ and $a_{ij{+}1}+1$,
respectively. Since $c_j(\Aq)<c_j(A)$, the minimality of $A$ implies
$x^\Aq\in X(\ell,d)$.  Now apply the polarization operator
$P_{jj{+}1}$ to $x^\Aq$. Then we get
$$\eqno{}
P_{jj{+}1}(x^\Aq)=(a_{ij{+}1}+1)x^A+a_{i{+}1j{+}1}x^{A_{i+1}}+\ldots+a_{\ell
j{+}1}x^{A_\ell}\in X(\ell,d).
$$
Here, $A_\nu$ is the matrix
$$\eqno{}
\pmatrix{
&\vdots&\vdots&\cr
\cdots&a_{ij}-1&a_{ij{+}1}+1&\cdots\cr
&\vdots&\vdots&\cr
\cdots&a_{\nu j}+1&a_{\nu j{+}1}-1&\cdots\cr
&\vdots&\vdots&\cr
}
$$
Since all $A_\nu$ are smaller than $A$ we obtain\footnote{In
characteristic 0, we are done at this point: $A$ can't exist.}
$a_{ij{+}1}+1=0$ in $k$. In particular, we get
$$\eqno{4}
a_{ij}\ge1\quad\hbox{and}\quad a_{ij{+}1}\ge p-1.
$$

Let $\ell':=\ell+1-i\ge j+1$. Then $a_{i\ell'}$ is on the diagonal
starting from $a_{1\ell}$:
$$\eqno{}
A=\pmatrix{
0        &\cdots&\cdots     &\cdots         &\cdots&0            &\cdots&a_{1\ell}\cr
\vdots   &      &&&      &\vdots       &\Ddots&\vdots\cr
0        &\cdots&a_{ij}&a_{ij{+}1}&\cdots&a_{i\ell'}   &      &\vdots\cr
\vdots   &      &      &          &\DDdots&\vdots       &      &\vdots\cr
\vdots   &      &      &\DDdots    &      &\vdots              &      &\vdots\cr
\vdots   &      &\DDdots&          &      &\vdots              &      &\vdots\cr
a_{\ell0}&a_{\ell1}&\cdots      &\cdots          &\cdots      &a_{\ell\ell'}&\cdots&a_{\ell\ell}\cr
}
$$
Let $A'$ be the submatrix of $A$
consisting of columns $0,\ldots,\ell'$. From \cite{E4} we get $r_i\ge p$
which implies $\ell'<\ell$ by \cite{E2}. Moreover we have $p-1\le
c_{j+1}\le c_{\ell'+1}\le\ldots\le c_\ell$. Hence
$$\eqno{3}
\|deg|A'=c_0+\ldots+c_{\ell'}\le(p-1)\ell-(p-1)(\ell-\ell')=(p-1)\ell'.
$$
Now observe that the first $i-1=\ell-\ell'$ rows of $A'$ are zero. Hence,
the monomial $x^{A'}$ ``lives'' on $(V')\^{\ell'+1}$ where
$\|dim|V'=\ell'$. Thus we can use induction and conclude that
$x^{A'}$ can be obtained by polarization from $(V')\^{\ell'}$. But
then the same polarization process produces $x^A$ from $V\^{\ell'}$
in contradiction to the choice of $A$.\qed

\Remarks: 1. The proof shows that under the given conditions it suffices
to apply just polarization operators and column permutations.

2. One can show that the monomial $x_{10}(x_{11}\ldots
x_{1\ell})^{p-1}$ cannot be obtained from polarization of polynomials
on fewer than $\ell+1$ copies of $V$. This shows that the given degree bound
is optimal.

\beginsection WeylTheorem. Weyl's theorem

Let $k$ and $p$ be as in the preceding section. We are going to apply
polarizations to invariant theory. Let $U$ and $V$ be two finite
dimensional representations of a finite group $G$. Then the
$GL(n)$\_action on $U\oplus V\^n$ commutes with $G$, hence we get an
$GL(n)$\_action and a notion of polarization on $k[U\oplus
V\^n]^G$. If $m\le n$ we can restrict invariants from $U\oplus V\^n$
to $U\oplus V\^m$. This process is kind of inverse to polarization and
is called {\it restitution}.

\Theorem Extend. Let $S$ be a generating set of $G$\_invariants on
$U\oplus V\^n$ where $n\ge\|max|(\|dim|V,{\beta_k(G)\over p-1})$. Then
the polarization of $S$ generates the ring of invariants on $U\oplus
V\^m$ for any $m\ge n$.

\Proof: First, let us remark that $G$ is linearly reductive. Indeed,
otherwise $\beta_k(G)=\infty$ and $m$ would not exist. 

Let $\cP_m\subseteq\|End|k[U\oplus V\^m]$ be the subalgebra generated
by all polarization and permutation operators. Let $A(m)$ be the
algebra generated by $\cP_m\cdot S\subseteq k[U\oplus V\^m]$. Since
$\cP_m$ commutes with $G$ we have $A(m)\subseteq k[U\oplus
V\^m]^G$. We have to show equality.

The formulas $P(fg)=P(f)g+fP(g)$ for any polarization operator and
$\pi(fg)=\pi(f)\pi(g)$ for any permutation imply that
$$\eqno{9}
\cP_m(fg)\subseteq\cP_m(f)\cP_m(g)
$$
for all $f,g$. This shows that $A(m)$ is a
$\cP_m$\_module. Let $d\le\beta_k(G)$. Then the assumption on $n$ and
\cite{Polar} imply that the map
$$\eqno{}
\cP_m\otimes k[U\oplus V\^n]_d\pfeil k[U\oplus V\^m]_d
$$
is surjective. By linear reductivity of $G$ we obtain
$$\eqno{}
\cP_m\cdot k[U\oplus V\^n]_d^G=k[U\oplus V\^m]_d^G.
$$
On the other hand,
$$\eqno{}
\cP_m\cdot k[U\oplus V\^n]_d^G=\cP_m\cdot k[S]_d\subseteq\cP_mA(m)_d=A(m)_d.
$$
Thus, we have proved that $A(m)$ contains all $G$\_invariants of
degree $\le\beta_k(G)$ which implies that it consists of all
invariants.\qed

\noindent This theorem provides a means to construct a representation
of $G$ which has the ``most general'' ring of invariants. To make this
more precise, we extend our notion of polarization. Let $U$ be a
$G$\_module. For any subset $S$ of $k[U]$ let $\|Pol|^U(S)$ be the
$\|Aut|^G(U)$\_module generated by it. Observe, that if $|G|\in
k^\Times$ then $\|Aut|^G(U)$ is just a product of general linear
groups corresponding to the isotypic components of $U$. There are two
more operations: If $V\subseteq U$ then $\|Res|_V^U(S)$ be the set of
restrictions of elements of $S$ to $V$. Conversely, if $U$ is a
quotient of $V$ let $\|Ext|^V_U(S)$ be the set of pullbacks. Now, we
define

\Definition: A $G$\_module $U$ has {\it universal invariants\/} if for
every finite dimensional $G$\_module $V$ the set
$$\eqno{}
\|Res|_V\|Pol|^{U\oplus V}\|Ext|_U^{U\oplus V}(S)
$$
generates $k[V]^G$ where $S$ is any generating set of $k[U]^G$.
\medskip\noindent
In other words, if $U$ has universal invariants then one obtains a
generating set of invariants for any other $G$\_module by the
following process: 1. Start with generators for $k[U]^G$. 2. Think of
them as invariants on $U\oplus V$. 3. Apply $\|Aut|^G(U\oplus V)$ to
$S$. 4. Restrict the ensuing invariants to $V$.

One of the main properties of a module with universal invariants is:

\Theorem. Assume $U$ has universal invariants. Then
$\beta_k(G)=\beta(k[U]^G)$.

\Proof: First, $\beta_k(G)$ is the supremum of the $\beta(k[V]^G$
where $V$ runs through all finite dimensional $G$\_modules
(\cite{vector}). Then the assertion follows from the fact that the
process of extension, polarization, and restriction does not change
degrees.\qed

\noindent Our main criterion for having universal invariants is:

\Theorem UniCrit. Assume $|G|\in k^\Times$. Let $U$ be a
representation of $G$ in which every irreducible module $M$ appears
with multiplicity at least $\|max|(\|dim|M,{\beta_k(G)\over
p-1})$. Then $U$ has universal invariants.

\Proof: Let $S$ be a generating set of $k[U]^G$ and let $V$ be a
finite dimensional $G$\_module. We show by induction on the number of
isotypic components of $V$ the more general statement that
$\|Aut|^G(U\oplus V)\cdot S$ generates $k[U\oplus V]^G$. We start with
$V=0$ where the assertion is trivial.

Assume that the simple module $M$ appears in $V$. Write $U=U'\oplus
M^{\oplus m}$ and $V=V'\oplus M^{\oplus n}$ such that $M$ is not
contained in $U'$ and $V'$. The set
$S':=\|Aut|^G(U\oplus V')\cdot S$ generates $k[U\oplus V']^G$ by the
induction hypothesis. We have $U\oplus V'=(U'\oplus V')\oplus
M^{\oplus m}$ while $U\oplus V=(U'\oplus V')\oplus M^{\oplus
m+n}$. \cite{Extend} implies that $k[U\oplus V]^G$ is generated by
$Gl(m+n)\cdot S'=Gl(m+n)\cdot \|Aut|^G(U\oplus V')\cdot
S=\|Aut|^G(U\oplus V)\cdot S$.\qed 

\Corollary. Assume $|G|\in k^\Times$. Then $G$ has a finite
dimensional representation with universal invariants. More precisely,
assume that $U$ contains every irreducible module with multiplicity at
least $|G|$. Then $U$ has universal invariants.

\Proof: Let $M$ be an irreducible representation of $G$. Since $M$ is a
quotient of the regular representation, we have $\|dim|M\le|G|$. On
the other hand, ${\beta_k(G)\over p-1}\le|G|$ by \cite{CB}.\qed

\Corollary regrep1. Assume $p>\beta_k(G)$. Then
the regular representation of $G$ has universal invariants.

\Proof: We may assume that $k$ is algebraically closed. Then the regular
representation contains each simple module $W$ with multiplicity
$\|dim|W$. The assertion follows from ${\beta_k(G)\over
p-1}\le1\le\|dim|W$.\qed

\beginsection improvement. A further improvement

One can try to obtain better results by generalizing the polarization
process. The point is that in positive characteristic, the commutator
algebra of $GL(V)$ in $\|End|(k[V^{\oplus m}]$ is not generated by the
image of $GL(m)$ (not even in a topological sense).

\def\cPq{{\overline\cP}}

In general, one can replace $\cP_m$ by any subalgebra $\cPq_m$ of endomorphisms
of $k[V\^m]$ which commutes with the $GL(V)$\_action and which
satisfies the multiplicativity property \cite{E9}. We will not pursue
this direction in any detail but illustrate it by the following easy fact:

\Theorem onedim. Assume $\|dim|V=1$ and
$\cPq_m:=\|End|^{GL(V)}(k[V\^m])$. Then for any $d\ge0$ holds
$$\eqno{}
k[V\^m]_d=\<k[V]_d\>_{\cPq_m}.
$$

\Proof: The group $GL(V)\cong GL(1)$ acts on $W_d:=k[V\^m]_d$ by the
character $t\mapsto t^d$. Hence,
$\cPq_m=\prod_{d=0}^\infty\|End|_kW_d$. This implies that any
non\_zero element generates $W_d$ as a $\cPq_m$\_module.\qed

\Remark: This theorem also serves as an example that $\cPq_m$ may be
bigger than $\cP_m$. In fact, $k[V\^m]_d$ is a simple $\cPq_m$\_module
but, in general, it is not a simple $GL(m)$\_module.
\medskip

\Lemma. Let $\cPq_m$ be as in \cite{onedim}. Then for any $f,g\in
k[V\^m]$ holds $\cPq_m(fg)\subseteq\cPq_m(f)\cPq_m(g)$.

\Proof: Let $f=\sum_d f_d$ be the decomposition into homogeneous
components. Since $\cPq_m$ contains the projection onto $k[V\^m]_d$ we
have $f_d\in\cPq_m(f)$. Hence, we may assume that $f$ and $g$ are
homogeneous of degree, say, $d$ and $e$, respectively. But then
$\cPq_m(f)=k[V\^m]_d$ and $\cPq_m(g)=k[V\^m]_e$, hence
$\cPq_m(f)\cPq_m(g)=k[V\^m]_{d+e}\supseteq\cPq_m(fg)$.\qed

Now assume $|G|\in k^\Times$. For an arbitrary $G$\_module $V$ we
generalize the polarization process as follows. Let
$V=\oplus_iM_i^{m_i}$ be the isotypic decomposition of $V$. Let
$\cPq_V\subseteq\|End|k[V]$ be the subalgebra which is generated by
$GL(m_i)$ if $\|dim|M_i>1$ and by $\cPq_{m_i}$ if $\|dim|M_i=1$. For a
subset $S\subseteq k[V]$ let $\overline{\rm Pol}^V(S):=\cPq_V(S)$.

We say a module $U$ has {\it weakly universal
invariants} if $k[V]^G$ is generated by
$$\eqno{}
\|Res|_V\overline{\rm Pol}^{U\oplus V}\|Ext|_U^{U\oplus V}(S)
$$
where $V$ is any finite
dimensional $G$\_module and $S$ is any generating set of
$k[U]^G$. Then we get the following analogue of \cite{UniCrit}

\Theorem wUniCrit. Assume $|G|\in k^\Times$Let $U$ be a representation
of $G$ in which every one\_dimensional module appears at least once
and every other irreducible module $M$ appears with multiplicity at
least $\|max|(\|dim|M,{\beta_k(G)\over p-1})$. Then $U$ has weakly
universal invariants.

\noindent The analogue of \cite{regrep1} is

\Corollary regu. Assume $k$ is algebraically closed. Let
$$\eqno{}
\ell:=\|inf|\{\|dim|W\mid W\hbox{ simple $G$\_module, }\|dim|W>1\}
\in[2,\infty].
$$
Assume $p\ge{\beta_k(G)\over\ell}+1$. Then the regular representation
$R$ of $G$ has weakly universal invariants. In particular,
$\beta_k(G)=\beta(k[R]^G)$.

\noindent If we combine this result with the bound of Domokos\_Heged\H
us\_Sezer we obtain the following strengthening of theorems of
Schmid, \cite{Sch}, and Smith, \cite{Sm}.

\Corollary regrep2. Assume $|G|\in k^\Times$ and $p\ge{3\over8}|G|+1$.
Then $\beta_k(G)=\beta(k[R]^G)$.

\Proof: If $G$ is cyclic then $\ell=\infty$. Otherwise $\ell\ge2$ and
$\beta_k(G)\le{3\over4}|G|$ by \cite{DomHegSez}.\qed

\Remark: The condition in \cite{Sch} is $\|char|k=0$. In \cite{Sm} it
is $p>|G|$.

\beginrefs

\L|Abk:DH|Sig:DH|Au:Domokos, M.; Heged\H us, P.|Tit:Noether's bound
for polynomial invariants of finite
groups|Zs:Arch. Math.|Bd:74|S:161--167|J:2000|xxx:-||

\L|Abk:Fl|Sig:Fl1|Au:Fleischmann, P.|Tit:The Noether bound in invariant
theory of finite groups|Zs:Adv. Math.|Bd:156|S:2000|J:23--32|xxx:-||

\L|Abk:Fleisch|Sig:Fl2|Au:Fleischmann, P.|Tit:On invariant theory of
finite groups|Zs:Preprint|Bd:-|S:49 pages|J:2002|xxx:-||

\L|Abk:Fo|Sig:Fo|Au:Fogarty, J.|Tit:On Noether's bound for polynomial
invariants of a finite group|Zs:Electron. Res. Announc. Amer. Math.
Soc.|Bd:7|S:5--7 (electronic)|J:2001|xxx:-||

\L|Abk:Noe|Sig:N|Au:Noether, E.|Tit:Der Endlichkeitssatz der Invarianten
endlicher Gruppen|Zs:Math. Ann.|Bd:77|S:89--92|J:1916|xxx:-||

\L|Abk:Ri|Sig:R|Au:Richman, D.|Tit:Invariants of finite groups over
fields of characteristic $p$|Zs:Adv. Math.|Bd:124|S:25--48|J:1996|xxx:-||

\Pr|Abk:Sch|Sig:Sch|Au:Schmid, B.|Artikel:Finite groups and invariant
theory|Titel:Topics in invariant theory (Paris, 1989/1990)|Hgr:M.-P.
Malliavin ed.|Reihe:Lecture Notes in
Mathematics|Bd:1478|Verlag:Springer\_Verlag|Ort:Berlin|S:35--66|J:1991|xxx:-||

\L|Abk:Se|Sig:Se|Au:Sezer, M.|Tit:Sharpening the generalized Noether
bound in the invariant theory of finite groups|Zs:J.
Algebra|Bd:254|S:252--263|J:2002|xxx:-||

\L|Abk:Sm|Sig:Sm|Au:Smith, L.|Tit:On a theorem of Barbara
Schmid|Zs:Proc. Amer. Math. Soc.|Bd:128|S:2199-2201|J:2000|xxx:-||

\B|Abk:Weyl|Sig:We|Au:Weyl, Hermann|Tit:The Classical Groups. Their
Invariants and Representations|Reihe:Princeton Mathematical
Series~{\bf 1}|Verlag:Princeton University
Press|Ort:Princeton|J:1939|xxx:-||

\endrefs

\bye